\documentclass[11pt,
]{article}                          

\usepackage{amssymb,amsfonts,amstext,amsmath,amsthm,latexsym,mathrsfs}
\usepackage{mathbbol}

\usepackage{mparhack}
\usepackage{makeidx}

\makeindex

\makeatletter
\if@twoside%
\else\fi%
\makeatother  

\input{ds.sty}

\begin{document}

\title
{OD elements of countable OD sets in the Solovay model
}

\author 
{
Vladimir~Kanovei\thanks{IITP RAS and MIIT,
  Moscow, Russia, \ {\tt kanovei@googlemail.com} --- contact author. 
}  
}

\date 
{\today}

\maketitle

\begin{abstract}
It is true in the Solovay model that every countable 
ordinal-definable set of sets of reals 
contains only ordinal-definable elements.
\end{abstract}

\parf{Introduction}
\las{cha1}

It is known that the existence of a non-empty $\od$ 
(ordinal-definable) set of reals $X$ with no \od\  
element is consistent with $\ZFC$; the set of all 
non-constructible reals gives a transparent example
in many generic models.
\bce
\rit{Can such a set\/ $X$ be countable}? 
\ece
This question was initiated and 
discussed at the \rit{Mathoverflow} website\snos
{\label{snos1}
\rit{Mathoverflow}, March 09, 2010. 
{\tt http://mathoverflow.net/questions/17608}. 
}
and at FOM\snos
{\label{snos2}%
FOM Jul 23, 2010.
{\tt http://cs.nyu.edu/pipermail/fom/2010-July/014944.html}}
.
In particular Ali Enayat (Footnote~\ref{snos2}) conjectured that 
the problem can be solved by the   
finite-support countable product $\plo$ (see \cite{ena})
of the Jensen ``minimal $\ip12$ real 
singleton forcing'' $\dP$ defined in \cite{jenmin} 
(see also Section 28A of \cite{jechmill}). 
We proved in \cite{kl:cds} 
that indeed, in a \dd\plo generic extension of $\rL$, 
the set of all reals \dd\dP generic over $\rL$ 
is a countable $\ip12$ set with no OD elements. 
Moreover there is a 
modification $\dP'$ of $\dP$ such that
it is true in a \dd{\dP'}generic extension of $\rL$ that
there is a $\ip12$ \dd\Eo equivalence class containing no
\od\ reals, \cite{kl:dec}.

On the other hand, one may ask 
do countable non-empty \od\ sets without 
\od\ elements exist in such a more typical generic extension  
as the Solovay model? 
We partially answer this question in the negative.

\bte
\lam{mt}
It is true in the Solovay model 
that every non-empty \od\ countable or finite set\/ 
$\cX$ {\bf of sets of reals}
necessarily contains an\/ \od\ element, and hence, in fact, 
consists of\/ $\od$ elements.
\ete

The Solovay model here is a model of $\zfc$ defined in 
\cite{solo} in which all projective 
(and generally all \rod, real-ordinal definable) 
sets of reals are Lebesgue measurable.
The case, when $\cX$ is a (non-empty \od\ countable) 
{\bf set of reals} in this theorem, is well 
known and is implicitly contained in the proof of the perfect 
set property by Solovay~\cite{solo}.
Hovever the proofs known for this particular case 
of sets of reals (as, \eg, in \cite{stern} or \cite{ksol}) 
do not work even for sets $\cX\sq\pws{\bn}$ 
(as in the theorem). 
In this paper, we present the proof of Theorem~\ref{mt}.

%

\parf{Notation} 
\las{tre}

We consider the constructible universe $\rL$ as the ground 
model by default. 
Suppose that $\Om$ is an inaccessible cardinal. 

\bbla
\lam{gset}
By a \rit{generic set} we'll always mean a \rit{filter}, 
that is, both pairwise compatible in itself and containing 
all weaker conditions.
\ebla

\bdf
\lam{lsf}
We represent the \rit{Levy -- Solovay forcing} associated 
with $\Om$ is the set $\ls$ of all partial maps 
\index{zzLS@$\ls$}
$p:\dom p\to\Om$ 
such that $\dom p\sq\Om\ti\om$ is a finite set and 
$p(\al,n)<\al$ whenever $\ang{\al,n}\in\dom p$.
Let $\abs p=\ens{\al}{\sus n\,(\ang{\al,n}\in\dom p)}$.
\index{zzpII@$\abs p$}

If $\ga\le\Om$ then $\ls_\ga=\ens{p\in\ls}{\abs p\sq\ga}$; 
\index{zzLSga@$\ls_\ga$}
in particular $\ls_\Om=\ls$.

If $p\in\ls$ and $\al<\Om$ then the \dd\al\rit{component} 
$\kom p\al$ of $p$ is a map defined on the set
$\dom{\kom p\al}=\ens{n}{\ang{\al,n}\in\dom p}\sq\om$ by 
$\kom p\al(n)=p(\al,n)$.
\edf

If $G\sq\ls$ 
is an \dd\ls generic set over $\rL$ then $\rL[G]$ is 
{\ubf the Solovay model}, to which Theorem~\ref{mt} refers.
\index{Solovay model}
The next lemma will be important below.

\ble
[reduction to ROD]
\lam{mtl}
It is true in the Solovay model that if\/ $\cX$ 
is a non-empty \od\ countable set and\/ $X\in\cX$ is\/ $\rod$ 
then\/ $X$ is\/ $\od$.
\ele

Thus somewhat surprisingly, 
it turns out that 
it suffices to prove the existence of a \rod\ 
(real-ordinal definable) element $X\in\cX$ 
in Theorem~\ref{mt}. 

\bpf
Arguing in the Solovay model, assume 
that 
$$
X=X_{p_0}=\ens{x}{\vpi(x,p_0)}\,, 
$$
where $\vpi$ is a formula with a real parameter $p_0\in\bn$ 
and hidden ordinal parameters.
The set $P=\ens{p\in\bn}{X_p\in\cX}$ is \od\ and contains 
$p_0$, and the equivalence relation, 
$p\rE q$ iff $X_p=X_q$ on $P$, is \od\ as well, and $\rE$ 
has at most countably many equivalence classes in $P$.
However it is known that, in the Solovay model, if an \od\ 
equivalence relation on $\bn$ has at most countably many 
equivalence classes then all its equivalence classes are 
\od, \cite{ksol,stern}. 
In particular $[p_0]_{\rE}$ is \od, and hence 
the set 
$X=X_{p_0}=\ens{x}
{\sus p\in[p_0]_{\rE}\,\vpi(x,p_0)}$
is \od.
\epf

\bdf
[ramified names]
\lam{nam1}
We'll use the ordinary ramified system of 
\dd\ls\rit{names} 
for differens sets in $\rL[G]$, so that 
$\bint UG$ will be the \dd G\rit{interpretation} of a 
\index{zzUGIIII@$\bint UG$}
name $U$ (basically, any set) defined by \dd\in rank 
induction by 
$$
\bint UG=\ens{\bint uG}{\sus p\in G\,(\ang{p,u}\in U)}\,.
$$
Then, if $G\sq\ls$ is generic over $\rL$ then
$\rL[G]=\ens{\bint UG}{U\in \rL}$.
\edf

Each set $x\in\rL$ has a canonical \dd\ls name 
$\na x\in\rL$, such that 
$\nt{\na x}G=x$ for any generic set $G\sq\ls$. 
Yet following common practice we shall identify $\na x$ 
with $x$ itself whenever possible.      

\bdf
[simple names]
\lam{nam2}
To somewhat simplify notation, we'll make 
use of a simpler system of names particularly for subsets 
of $\ls$. 
Let $\lsn=\pws{\ls\ti\ls}$, and if 
\index{zzN@$\lsn$}
$t\in\lsn$ and $G\sq\ls$ then 
$\nt tG=\ens{q}{\sus p\in G\,(\ang{p,q}\in t)}\sq\ls$.
\index{zztGII@$\nt tG$}

Thus $\lsn$ consists of all \dd\ls names for subsets of $\ls$.

If $\ga<\Om$ then let 
$\lsn_\ga=\pws{(\ls_\ga)\ti(\ls_\ga)}$, so that any 
\index{zzNga@$\lsn_\ga$}
$t\in\lsn_\ga$ is a \dd{\ls_\ga}name for a subset of 
$\ls_\ga$.
\edf 

The name $\uG=\ens{\ang{p,p}}{p\in\ls}$ belongs to $\lsn$, 
\index{zzGund@$\uG$}
and $\uG[G]=G$.

\parf{Double names}
\las{qua}   

In many cases below, we'll consider pairs of 
\dd\ls generic sets $G,G'\sq\ls$ over $\rL$, 
such that $\rL[G]=\rL[G']$; note that this is 
not a \dd{(\ls\ti\ls)}generic pair! 
Similar pairs will be considered for the forcing 
notions $\ls_\ga$ ($\ga<\Om$) instead of $\ls$. 
The next definition introduces a useful tool related 
to such pairs.

\bdf
\lam{dq}
In $\rL$, if $\ga\le\Om$ then any pair 
$a=\ang{\tl{a},\tr{a}}$ of names 
\index{zztalef@$\tl a$}%
\index{zztarig@$\tr a$}%
$\tl{a},\tr{a}\in\lsn_\ga$ will be called a 
\rit{double-name}. 
\index{double-name}
Let $\trn_\ga$ consist of all 
\index{zzDNga@$\trn_\ga$}
double-names $a=\ang{\tl{a},\tr{a}}$ such that  
$\tl a\ne\pu$, $\tr a\ne\pu$, and
\ben
\nenu
\itla{dq1} 
if $p\in\dom\tl a$ then $p$ \dd{\ls_\ga}forces: \ 
\begin{minipage}[t]{\textwidth}
(a)
$\tl{a}[\uG]$ is 
\dd{\ls_\ga}generic, and \\[1ex] 
(b)
$\uG=\tr{a}[\tl{a}[\uG]]$;\vom
\end{minipage}

\itla{dq2}
if $p\in\dom\tr a$ then $p$ \dd{\ls_\ga}forces: \  
\begin{minipage}[t]{\textwidth}
(a)
$\tr{a}[\uG]$ is 
\dd{\ls_\ga}generic, and \\[1ex] 
(b)
$\uG=\tl{a}[\tr{a}[\uG]]$.
\end{minipage}
\een 
Define $\trn=\bigcup_{\ga<\Om}\trn_\ga$; 
\index{zzDN@$\trn$}
this is different from $\trn_\Om$. 
It follows from \ref{dq1} or \ref{dq2} that for any 
$a\in\trn$ there is a \rit{unique} $\ga=\abt{a}<\Om$ such that 
\index{zzaII@$\abt{a}$}
$a\in\trn_\ga$.
\edf  

Note that all sets $\lsn_\ga$ and $\trn_\ga$ belong to $\rL$.

\vyk{
\ble
[\bf do I need it?]
\lam{GGc}
If\/ $\ga\le\Om$ and\/ $p,q\in\ls_\ga$ 
then there is a double-name\/ ${a}\in\trn_\ga$ such that\/ 
$p\in\dom\tl a$ and\/ $q\in\dom\tr a$. 
\ele
\bpf
If $\al<\ga$ then $N^p_\al=\ens{n}{\ang{\al,n}\in\dom p}$ 
and $N^q_\al=\ens{n}{\ang{\al,n}\in\dom q}$ are finite 
sets; let $h_\al$ be the order-increasing bijection of 
$\om\bez N^p_\al$ onto $\om\bez N^q_\al$.
Note that $N^p_\al=N^q_\al=\pu$ and $h_\al$ is the identity 
on $\om$ for all but finite $\al$. 
Define $h(\al,n)=\ang{\al,h_\al(n)}$; this is a bijection 
of the cofinite set $D_p=(\ga\ti\om)\bez \dom p$ onto 
the cofinite   
$D_q=(\ga\ti\om)\bez \dom q$, equal to 
the identity on $(\ga\bez(\abs p\cup\abs q))\ti\om$.

If $r\in\ls_\ga$ and $\dom r\sq D_p$ then let 
$r^h\in\ls_\ga$ be defined by 
$\dom {r^h}=\ens{h(\al,n)}{\ang{\al,n}\in\dom r}$ 
and $r^h(h(\al,n))=r(\al,n)$ for all $\ang{\al,n}\in\dom r$.
Finally let 
$\tl a=\ens{\ang{p\cup r,q\cup r^h}}
{r\in\ls_\ga\land \dom r\sq D_p}$
and $\tr a=(\tl a)^{-1}$.
\epf
}

\ble
\lam{io}
Assume that\/ $\ga\le\Om$ and\/ ${a}\in\trn_\ga$.
Then$:$
\ben
\renu
\itla{io1}
if\/ $\Gl\sq\ls_\ga$ 
is an\/ \dd{\ls_\ga}generic set and\/ 
$\Gl\cap\dom {\tl a}\ne\pu$ then 
$\Gr=\tl{a}[\Gl]$ is\/ \dd{\ls_\ga}generic, 
$\Gr\cap\dom{\tr a}\ne\pu$, and\/ $\Gl=\tr{a}[\Gr]\,;$

\itla{io2}
if\/ $\Gr\sq\ls_\ga$ 
is\/ \dd{\ls_\ga}generic and\/ 
$\Gr\cap\dom{\tr a}\ne\pu$ then 
$\Gl=\tr{a}[\Gr]$ is\/ \dd{\ls_\ga}generic, 
$\Gl\cap\dom {\tl a}\ne\pu$, $\Gr=\tl{a}[\Gl]\,.$
\qed
\een
\ele

Thus each 
${a}\in\trn_\ga$ induces a bijection between 
all \dd{\ls_\ga}generic sets $G\sq\ls_\ga$ satisfying 
$G\cap\dom{\tl a}\ne\pu$ and those satisfying 
$G\cap\dom{\tr a}\ne\pu$.

\bcor
\lam{qq'}
If\/ $\ga\le\Om$, ${a}\in\trn_\ga$,
$\ang{q,p}\in\tr a$, and\/ $q\sq q'\in\ls_\ga$
then there is a condition\/ $p'\in\ls_\ga$ 
compatible with\/ $p$ and such that\/ 
$\ang{p',q'}\in\tl a$.
\ecor
\bpf
Let $\Gr\sq\ls_\ga$ be a generic set containing $q'$, 
hence containing $q$ as well.
Then $\Gl=\tr{a}[\Gr]$ is a \dd{\ls_\ga}generic set 
containing $p$, 
and $\Gr=\tl{a}[\Gl]$ by Lemma{io}. 
As $q'\in\Gr$, there is a condition $p'\in\Gl$ 
such that $\ang{p',q'}\in\tl a$. 
As $p$ also belongs to $\Gl$, $p,p'$ are compatible.
\epf

\parf{Full, regular, equivalent names}
\las{fre}   

Recall that a set $D\sq\ls_\ga$ is \rit{dense} 
\index{set!dense}
if for any $p\in\ls_\ga$ there is $q\in D$ 
with $p\sq q$, and is \rit{open} if 
\index{set!open}
$(p\in D\land p\sq q\in\ls_\ga)\imp q\in D$. 

\bdf
\lam{fn}
Let $\ga\le\Om$. 
A name $t\in\lsn_\ga$ is \rit{full} if the set $\dom t$ 
is dense in $\ls_\ga$. 
A double-name $a\in\trn_\ga$ is \rit{full} if such are 
\index{name!full}%
\index{double-name!full}%
the names $\tl a$ and $\tr a$.

A name $t\in\lsn_\ga$ is \rit{regular}, 
if the following holds: 
if $p,q\in\ls_\ga$ and $p$ \dd{\ls_\ga}forces $q\in t[\uG]$ 
then $\ang{p,q}\in t$. 
\index{name!regular}%
\index{double-name!regular}%
In particular, in this case, if $\ang{p,q}\in t$ and 
$p\sq p'\in\ls_\ga$ then $\ang{p',q}\in t$, too. 
A double-name $a\in\trn_\ga$ is \rit{regular},
if so are both components $\tl a$ and $\tr a$.
Define the \rit{regular hull}  
\index{regular hull}%
\index{zzrht@$\rh t$}%
$$
\rh t=\ens{\ang{p,q}\in\ls_\ga\ti\ls_\ga}
{\text{$p$ \dd{\ls_\ga}forces $q\in t[\uG]$}}. 
\index{zztrh@$\rh t$} 
$$
of any $t\in\lsn_\ga$. 
If $a\in\trn_\ga$ then let 
$\rh a=\ang{\rh{\tl a},\rh{\tr a}}$.
\edf

\ble
\lam{fnd}
Assume that $\ga\le\Om$ and\/ $a\in\trn_\ga$ is full.  
Then\/ $\ran{\tl a}=\ran{\tr a}=\ls_\ga$, and if\/ 
$G\sq\ls_\ga$ is\/ \dd{\ls_\ga}generic then so are 
$\tl a[G]$ and\/ $\tr a[G]$.
\ele
\bpf
To prove the genericity claim note that if say 
$\dom{\tl a}$ is dense then any generic set $G\sq\ls_\ga$ 
intersects $\dom{\tl a}$, then use Lemma~\ref{io}. 
To prove the first claim, let $q\in\ls_\ga$. 
Consider a generic set $\Gr\sq\ls_\ga$ containing $q$. 
Then $G\cap \dom{\tr a}\ne\pu$, see above. 
It follows that $\Gl=\tr a[\Gr]$ is generic and 
$\Gr=\tl a[\Gl]$ by Lemma~\ref{io}. 
But $q\in \Gr$, hence $q\in \ran{\tl c}$.
\epf

\bdf
\lam{eqn}
Names $s,t\in\lsn_\ga$ are \rit{equivalent} if 
$s[G]=t[G]$ for any generic set $G\sq\ls_\ga$, or 
equivalently, 
if any $p\in\ls_\ga$ \dd{\ls_\ga}forces $s[\uG]= t[\uG]$. 
\index{name!equivalent}%
\index{double-name!equivalent}%
Double-names $a,b\in\trn_\ga$ are 
\rit{equivalent} if $\tl b$, $\tr b$ are equivalent to 
resp.\ $\tl a$, $\tr a$.
\edf

\ble
\lam{rh}
Assume that\/ $\ga\le\Om$. 
Then$:$
\ben
\renu
\itla{rh1}
if\/ $t\in\lsn_\ga$ then\/ $\rh t$ is 
regular and equivalent to\/ $t\,;$  

\itla{rh2}
if\/ $a\in\trn_\ga$ then\/ $\rh a\in\trn_\ga$, 
$a\les\rh a$, and\/ 
$\rh a$ is equivalent to\/ $a$ --- therefore the set\/ 
$\trr_\ga=\ens{b\in\trn_\ga}{b\,\text{\rm\ is regular}}$ 
is dense in\/ $\trn_\ga\,;$ 

\itla{rh3}
if\/ $a,b\in\trn_\ga$ then\/ 
$a$ is equivalent to\/ $b$ iff\/ $\rh a=\rh b$. 
\een
\ele
\bpf
\ref{rh1} 
To establish the equivalence, assume that $G\sq \ls_\ga$  
is generic and $q\in\rh t[G]$. 
Then there is $p\in G$ such that $\ang{p,q}\in\rh t$. 
By definition 
$p$ \dd{\ls_\ga}forces $q\in t[\uG]$. 
But then $q\in t[G]$, as required.
To establish the regularity, assume that $p,q\in\ls_\ga$, 
and $p$ \dd{\ls_\ga}forces $q\in{\rh t}[\uG]$ 
--- therefore $p$ \dd{\ls_\ga}forces $q\in{t}[\uG]$ by 
the equivalence already proved. 
Then by definition $\ang{p,q}\in \rh t$. 

\ref{rh2} follows from \ref{rh1}. 
The direction $\mpi$ in \ref{rh3} immediately follows from 
\ref{rh2}. 
To prove the opposite direction, it suffices to show that 
if names $s,t\in\lsn_\ga$ are equivalent then $\rh s=\rh t$.
Assume that $\ang{p,q}\in\rh s$. 
By definition $p$  \dd{\ls_\ga}forces $q\in s[\uG]$. 
Then, as $s,t$ are equivalent, $p$ also forces $q\in t[\uG]$.
It follows that $\ang{p,q}\in\rh t$, as required.
\epf

\bpri
\lam{idg}
If $\ga<\Om$ then let 
$t_\ga=\ens{\ang{p,q}}{p,q\in\ls_\ga\land q\sq p}$ and 
$\ide\ga=\ang{t_\ga,t_\ga}$.
\index{zzidga@$\ide\ga$}%
\index{name!idga@$\ide\ga$}%
Then $\ide\ga\in\trn_\ga$ is a full regular double-name 
and $\tl{\ide\ga}[G]=\tl{\ide\ga}[G]=G$ for any 
\dd{\ls_\ga}generic set $G\sq\ls_\ga$: 
the \rit{identity\/} name. 
\epri


\parf{Double-name representation theorem}
\las{dnr}   

The next theorem shows that the double-name tool 
adequately represents the case of a pair of 
\dd\ls generic sets $G,G'\sq\ls$  
such that $\rL[G]=\rL[G']$.

\bte
\lam{GG}
Assume that\/ $\ga\le\Om$, $\Gl,\Gr\sq\ls_\ga$ 
are\/ \dd{\ls_\ga}generic sets over\/ $\rL$, 
and\/ $\rL[\Gl]=\rL[\Gr]$.
Then there is a full regular double-name\/ 
$c\in\trn_\ga$ such that\/ 
$\Gr=\tl{c}[\Gl]$, $\Gl=\tr{c}[\Gr]$, and\/ 
$\tl c=\tr c$.
\ete
\bpf
If $\Gl=\Gr$ then it suffices to define $c$ by 
$\tl c=\tr c=\ide\ga$. 
Therefore assume that $\Gl\ne\Gr$. 
Then there exist conditions $p\lef\in\Gl$ and $p\rig\in\Gr$ 
incompatible in $\ls_\ga$.
By a basic forcing theorem, there exist names 
$s\lef,s\rig\in\lsn_\ga$ such that 
$\Gr=s\lef[\Gl]$, $\Gl=s\rig[\Gr]$, and every condition 
$p\in\dom{s\lef}$ satisfies $p\lef\sq p$ while 
every condition $q\in\dom{s\rig}$ satisfies $p\rig\sq q$.
It is not true immediately that 
$\ang{s\lef,s\rig}\in\trn_\ga$; we need to somewhat 
modify the  names by shrinking.

We can wlog assume that $s\lef$ and $s\rig$ are 
regular; as otherwise we can replace them by resp.\ 
$\rh{s\lef}$ and $\rh{s\rig}$ and use Lemma~\ref{rh}\ref{rh1}. 

Define $a=\ang{\tl a,\tr a}$,  
where $\tl a$ consists of all pairs $\ang{p,q}\in s\lef$ 
such that 
\bce
$p$ \dd{\ls_\ga}forces that  
$s\lef[\uG]$ is \dd{\ls_\ga}generic and $\uG=s\rig[s\lef[\uG]]$,
\ece 
and $\tr a$ consists of all pairs $\ang{q,p}\in s\rig$ 
such that 
\bce
$q$ \dd{\ls_\ga}forces that $s\rig[\uG]$ is 
\dd{\ls_\ga}generic and $\uG=s\lef[s\rig[\uG]]$;
\ece 
then $\pu\ne\tl a\sq s\lef$ and $\pu\ne\tr a\sq s\rig$. 

We claim that $a\in\trn_\ga$, and still 
$\Gr=\tl a[\Gl]$ and $\Gl=\tr a[\Gr]$. 

\ble
\lam{GGle} 
If\/ $\Hl$ is an\/ \dd{\ls_\ga}generic set and 
$\Hl\cap\dom{\tl a}\ne\pu$ then\/ $\tl a[\Hl]=s\lef[\Hl]$.
Similarly if\/ $\Hr$ is an\/ \dd{\ls_\ga}generic set and 
$\Hr\cap\dom{\tr a}\ne\pu$ then\/ $\tr a[\Hr]=s\lef[\Hr]$.
\ele
\bpf[lemma]
By construction $\tl a[\Hl]\sq =s\lef[\Hl]$.
Consider any $q\in s\lef[\Hl]$, so that there is $p\in\Hl$ 
with $\ang{p,q}\in s\lef$.
On the other hand, as $\Hl\cap\dom{\tl a}\ne\pu$, there 
is a condition $p'\in\Hl$ with $p\sq p'$ which 
\dd{\ls_\ga}forces that $s\lef[\uG]$ is 
\dd{\ls_\ga}generic and $\uG=s\rig[s\lef[\uG]]$. 
Then $\ang{p',q}\in \tl a$ by the regularity assumption, 
and we have $q\in \tl a[\Hl]$. 
\epF{Lemma}

Now to check \ref{dq}\ref{dq1} for $a$ 
let $\Hl$ be an \dd{\ls_\ga}generic set and 
$\Hl\cap\dom{\tl a}\ne\pu$. 
Then $\tl a[\Hl]=s\lef[\Hl]$ by the lemma. 
Therefore $\Hr=\tl a[\Hl]$ is \dd{\ls_\ga}generic 
and $\Hl=s\rig[\Hr]$ by the definition of $\tl a$. 
Thus $s\rig[\Hr]$ is generic and 
$s\lef[s\rig[\Hr]]=\Hr$ by construction. 
This is forced by some $q\in\Hr$. 
On the other hand, as $\Hl=s\rig[\Hr]\ne\pu$, there 
exists some $q'\in \Hr\cap\dom s\rig$. 
We can assume that $q'\sq q$. 
Then $q\in\dom s\rig$, too, by the regularity 
assumption, and hence $q\in\dom \tr a$, 
and $\Hr\cap\dom \tr a\ne\pu$. 
We conclude that $\tr a[\Hr]=s\rig[\Hr]=\Hl$, by 
the lemma.
Finally $\tr a[\tl a[\Hl]]=\Hl$; 
this ends the verification of \ref{dq}\ref{dq1} 
for $a$.

Thus $a\in\trn_\ga$. 
In addition, by the choice of $s\lef$ and $s\rig$, some 
$p\in \Gl$ forces that 
``$s\lef[\uG]$ is generic and $\uG=s\rig[s\lef[G]]$''. 
Then $p\in\dom s\lef$, $p\in\dom \tl a$,
$\Gl\cap\dom \tl a\ne\pu$, and $\tl a[\Gl]=s\lef[\Gl]=\Gr$,
as above.  
Similarly we have $\Gl=\tr a[\Gr]$.

To fix the regularity condition of the theorem, 
let $b=\rh a$; then still 
$b\in\trn_\ga$, $\Gr=\tl b[\Gl]$, $\Gl=\tr b[\Gr]$, 
and $b$ is regular, by Lemma~\ref{rh}.

It is not necessarily true, of course, that sets 
$\dom{\tl b}$ and $\dom{\tr b}$ are dense.
To fix this shortcoming, we define 
$$
W=\ens{p\in\ls_\ga}
{\kaz q\in \dom{\tl b}\cup\dom{\tr b}\,
(p\:\text{ is incompatible with }\,q)}
$$
and let $c=\ang{\tl c,\tr c}$, where
$ 
\tl c=\tr c=\tl b\cup\tr b\cup
\ens{\ang{p,q}}{p\in W\land q\sq p}\,.
$ 

The set  
$\dom{\tl c}=\dom{\tr c}=\dom{\tl b}\cup\dom{\tr b}\cup W$ 
is dense in $\ls_\ga$ by construction.
We claim that $c\in\trn_\ga$. 
Indeed let $\Hl\sq\ls_\ga$ be an 
\dd{\ls_\ga}generic set. 
Then $\Hl\cap \dom{\tl c}\ne\pu$. 
But $\dom{\tl c}=\dom{\tl b}\cup\dom{\tr b}\cup W$.\vom

{\it Case 1\/}: 
$\Hl\cap \dom{\tl b}\ne\pu$.
Then $\Hl\cap \dom{\tr b}=\pu$ since if $p'\in\dom{\tl b}$ 
and $q'\in\dom{\tl b}$ then $p',q'$ are incompatible by 
the original choice of $p\lef,p\rig$.
We also have $\Hl\cap W=\pu$ by obvious reasons. 
It follows that $\tl c[\Hl]=\tl b[\Hl]$, and hence 
$\Hr=\tl c[\Hl]$ is an \dd{\ls_\ga}generic set and 
$\Hl=\tr b[\Hr]$, because $b\in\trn_\ga$. 
In particular $\Hr\cap\dom{\tl b}\ne\pu$, so that 
$\tr c[\Hr]=\tr b[\Hr]$, as above.\vom  

{\it Case 2\/}: 
$\Hl\cap \dom{\tr b}\ne\pu$, similar.\vom  

{\it Case 3\/}: 
$\Hl\cap W\ne\pu$.
Then $\Hl\cap \dom{\tl b}=\Hl\cap \dom{\tr b}=\pu$ as above. 
It follows that $\tl c[\Hl]=\tr c[\Hl]=\Hl$.\vom

Thus indeed $c\in\trn_\ga$, 
$\tl c=\tr c$,    
the set $\dom{\tl c}=\dom{\tr c}$ is open dense 
in\/ $\ls_\ga$,  
and the arguments above (Case 1) also imply that 
$\Gr=\tl{c}[\Gl]$, $\Gl=\tr{c}[\Gr]$.
Moreover, $c$ inherits the regularity of $b$.
\epf

\parf{Extensions}
\las{ext}   

\bdf
[extension]
\lam{d:e}                  
Suppose that $a,b$ are double-names. 
We say that 
${b}$ \rit{extends} ${a}$, 
\index{extends}%
in symbol ${a}\les{b}$, if just 
\index{zzzles@$\les$}%
$\tl{a}\sq\tl{b}$ and $\tr{a}\sq\tr{b}$.
\edf

\ble
[in $\rL$]
\lam{ext1}
If\/ $\ba<\ga\le\Om$ and\/ ${a}\in\trn_\ba$, then there is a 
double-name\/ ${b}\in\trn_\ga$ which
extends\/ ${a}$.
\ele     
\bpf
Let $\tl{b}$ consist of all pairs 
$\ang{p\cup r,q\cup r}$, where $\ang{p,q}\in \tl{a}$ and 
$r$ is a condition in $\ls_\ga$ satisfying 
$\abt r \sq\ga\bez\ba$;
let $\tr{b}$ be defined the same way.

This can be explained as follows. 
Suppose that $\Gl\sq\ls_\ga$ is a \dd{\ls_\ga}generic 
set containing $p\lef$. 
Then the factors 
$G'\lef=\Gl\cap\ls_\ba$ and $G''\lef=\Gl\cap\ls_{\ga\bez\ba}$ 
are resp.\ \dd{\ls_\ba}generic and \dd{\ls_{\ga\bez\ba}}generic, 
and $\Gl$ can be identified with $G'\lef\ti G''\lef$ by the 
product forcing theorem.
Then by definition the set $\Gr=\tl{b}[\Gl]$ has the form 
$G'\rig\ti G''\rig$, where $G'\rig=\tl{a}[G'\lef]$ while simply 
$G''\rig=G''\lef$. 
The genericity of $\Gr$ easily follows.
\epf

\bdf
[restriction]
\lam{d:r}
Let $\al<\ba\le \Om$. 
If $t\in\ls_\ba$ then define 
$t\res\al=t\cap {(\ls_\al\ti\ls_\al)}$; 
$t\res\al\in \lsn_\al$.
\index{zztalI@$t\res\al$}%
\index{zzaalI@$a\res\al$}%
\index{restriction!talI@$t\res\al$}%
\index{restriction!aalI@$a\res\al$}%
If ${a}\in\trn_\ba$, then let 
${a}\res\al=\ang{\tl{a}\res\al,\tr{a}\res\al}$. 
\edf

It is not asserted that \rit{always} ${a}\res\al\in\trn_\al$!

\ble
\lam{0ext}
If, in\/ $\rL$, 
$\al<\ba\le\Om$, ${a}\in\trn_\al$, ${b}\in\trn_\ba$, 
and\/ ${a}\les{b}$, then 
\ben
\renu
\itla{0ext1}
if\/ $\Gl\sq\ls_\ba$ is an\/ \dd{\ls_\ba}generic set 
then\/ 
{\rm(a)} 
$H\lef=\Gl\cap\ls_\al$ is\/ 
\dd{\ls_\al}generic, and\/ 
{\rm(b)} 
if\/ $\Hl\cap\dom{\tl a}\ne\pu$ 
then\/ $\tl{a}[H\lef]=\tl{b}[\Gl]\cap\ls_\al\,;$ 

\itla{0ext2}
if\/ $\Gr\sq\ls_\ba$ is an\/ \dd{\ls_\ba}generic set 
then\/ 
{\rm(a)} 
$H\rig=\Gr\cap\ls_\al$ is\/ 
\dd{\ls_\al}generic, and 
{\rm(b)} 
if\/ $\Gr\cap\dom\tr a\ne\pu$  
then\/ $\tr{a}[H\rig]=\tr{b}[\Gr]\cap\ls_\al\,;$ 

\itla{0ext3}
$c=b\res\al$ belongs to\/ $\trn_\al$ and\/ $a\les c\les b$.
\een
\ele     
\bpf

\ref{0ext1}(a)
That $H\lef$ is generic holds by the product forcing theorem. 

\ref{0ext1}(b)
If $\Hl\cap\dom{\tl a}\ne\pu$ then 
$\Gl\cap\dom{\tl b}\ne\pu$, and hence
the sets $\Gr=\tl{b}[\Gl]$ and $H\rig=\tl{a}[H\lef]$ 
are generic sets in resp.\ $\ls_\ba$ and $\ls_\al$ by 
Lemma~\ref{io}, and $H\rig\sq \Gr$ since ${a}\les{b}$.
Therefore $H\rig\sq H'\rig=\Gr\cap\ls_\al$. 
However $H'\rig$ is \dd{\ls_\al}generic
by the product forcing. 
Thus both $H\rig\sq H'\rig$ are generic sets, 
hence easily $H\rig=H'\rig$ as required. 

\ref{0ext3} 
To check \ref{dq}\ref{dq1}(a) for some $p\in\dom{\tl c}$, 
consider any \dd{\ls_\al}generic set $\Hl\sq\ls_\al$ 
containing $p$ and extend it to a 
\dd{\ls_\ba}generic set $\Gl\sq\ls_\al$ so that 
$\Hl=\Gl\cap\ls_\al$. 
The (generic by Lemma~\ref{io}) 
sets $\Hr=\tl a[\Hl]$ and $\Gr=\tl b[\Gl]$ satisfy 
$\Hr=\Gr\cap\ls\ba$ by \ref{0ext1}. 
On the other hand $\Hr\sq \tl c[\Hl]\sq\Gr\cap\ls_\ba$, 
hence $\tl c[\Hl]=\Hr$ is generic, as required. 
The verification of \ref{dq}\ref{dq1}(b) also is very 
simple.
\epf

\ble
\lam{erh}
In\/ $\rL$, assume that\/ $\al<\ba\le\Om$. 
Then$:$
\ben
\renu
\itla{erh1}
if\/ $s\in\lsn_\al$, $t\in\lsn_\ba$, 
and\/ ${s}\sq{t}$, then\/ $\rh s\sq\rh t\,;$ 

\itla{erh2}
therefore if\/ ${a}\in\trn_\al$, ${b}\in\trn_\ba$, 
and\/ ${a}\les{b}$, then\/ $\rh a\les\rh b\,;$ 

\itla{erh3}
if\/ ${b}\in\trn_\ba$ is regular and\/ 
$a={b}\res\al\in\trn_\al$ then\/ $a$ is regular, too.
\een
\ele
\bpf
\ref{erh1}
Suppose that $\ang{p',q}\in\rh s$, \ie, $p',q\in\ls_\al$ 
and there is a condition $p\sq p'$ which 
\dd{\ls_\al}forces that $q\in t[\uG]$. 
Prove that $p$ also \dd{\ls_\ba}forces $q\in t[\uG]$. 
Let a set $\Gl\sq\ls_\ba$ be a set 
\dd{\ls_\ba}generic over $\rL$ and containing $p$; 
prove that $q\in\Gr=t[\Gl]$.
The set $\Hl=\Gl\res\ls_\al$ is \dd{\ls_\al}generic
by Lemma~\ref{0ext} and still $p\in\Hl$, hence 
$q\in s[\Hl]\sq t[\Gl]=\Gr$, as required.

\ref{erh3}
Assume that $p,q,p'\in\ls_\al$, $p\sq p'$ and $p$ 
\dd{\ls_\al}forces $q\in \tl a[\uG]$; we have to prove 
that $\ang{p',q}\in \tl a$. 
As $a={b}\res\al$, 
it suffices to show that $\ang{p',q}\in \tl b$.
The same argument based on Lemma~\ref{0ext} shows that 
$p$ also \dd{\ls_\ba}forces $q\in \tl a[\uG]$. 
Therefore $\ang{p',q}\in \tl b$ since $b$ is regular.
\epf

\parf{Increasing sequences}
\las{inc}   

Suppose that a set $\Ga\sq\trn$ is pairwise \dd\les 
compatible. 
Then define the double-name  
$\xa=\bigvee\Ga$ by 
\index{zzGabigvee@$\bigvee\Ga$}%
\index{zzzbigveeGa@$\bigvee\Ga$}%
$\tl\xa=\bigcup_{{a}\in\Ga}\tl{a}$, \
$\tr\xa=\bigcup_{{a}\in\Ga}\tr{a}$. 

\ble
[in $\rL$]
\label{3ext}
\ben
\renu
\itla{3ext1}
If\/ $\la<\Om$ is a limit ordinal and\/ 
\imar{3ext}
$\sis{a_\xi}{\xi<\la}$ is a\/ \dd\les increasing sequence 
in\/ $\trn$ then\/ $A=\bigvee\ens{a_\xi}{\xi<\la}$ 
belongs to\/ $\trn\;;$\vom

\itla{3ext2}
therefore the set\/ $\trn=\bigcup_{\ga<\Om}\trn_\ga$ 
is\/ \dd\Om closed in the sense of\/ $\les\;;$\vom 

\itla{3ext3}
if\/ $\sis{a_\xi}{\xi<\Om}$ is a strictly\/ 
\dd\les increasing sequence 
in\/ $\trn$ then the double-name\/ 
$A=\bigvee\ens{a_\xi}{\xi<\la}$ 
belongs to\/ $\trn_\Om$.
\een
\ele  
\bpf
\ref{3ext1} 
Suppose that  
$\sis{\ga_\xi}{\xi<\la}$ is a strictly increasing sequence 
of ordinals $\ga_\xi<\Om$,  
and double-names 
${a}_\xi=\ang{t\lef^\xi,t\rig^\xi}\in\trn_{\ga_\xi}$ 
form a strictly \dd\les increasing sequence: if 
$\xi<\eta<\la$ then $\tl\xi\sq\tl\eta$ and $\tr\xi\sq\tr\eta$.
Let $\tl A=\bigcup_{\xi<\la}t\lef^\xi$,  
$\tr A=\bigcup_{\xi<\la}t\rig^\xi$, 
and $\ga=\sup_{\xi<\la}\ga_\xi$. 
We claim  that $A=\ang{\tl A,\tr A}\in\trn_{\ga}$.  

Let's verify \ref{dq}\ref{dq1}. 
Assume that $\Gl\sq\ls_\ga$ is a generic set 
containing some $p\in\dom{\tl A}$; 
we have to prove that $\Gr=\tl A[\Gl]$ 
is \dd{\ls_\ga}generic and $\Gl=\tr A[\Gr]$. 
Note first of all that each set  
$\Gl^\xi=\Gl\cap\ls_{\ga_\xi}$, $\xi<\la$, 
is \dd{\ls_{\ga_\xi}}generic by the product forcing theorem, 
and $p$ belongs to some $\dom{\tl {a_\za}}$, $\za<\Om$.
We can assume that $\za=0$ 
(otherwise simply cut all double-names $a_\xi$, $\xi<\za$). 
Then $p\in\dom{\tl {0}}$, 
therefore $p\in\dom{\tl {\xi}}$ for all $\xi<\Om$. 
It follows that each set 
$\Gr^\xi=\tl\xi[\Gl^\xi]\sq\ls_{\ga_\xi}$ is 
\dd{\ls_{\ga_\xi}}generic, 
$\Gr^\xi\cap \dom{\tr\xi}\ne\pu$,  
and $\Gl^\xi=\tr\xi[\Gr^\xi]$, 
by Lemma~\ref{io}.
And as $\Gr=\bigcup_{\xi<\la}\Gr^\xi$, 
we conclude that at least $\Gr$ is a filter in $\ls_\ga$ 
and $\Gl=\tr A[\Gr]$, that is, \ref{dq}\ref{dq1}(b). 

\vyk{
Further, the sets $\Gr^\xi=t\lef^\xi[\Gl^\xi]\sq\ls_{\ga_\xi}$ 
satisfy $\Gr=\bigcup_{\xi<\la}\Gr^\xi$. 
It follows that $p\rig$ belongs to some $\Gr^\za$, 
and then there is a condition $p_0\in\Gl^\za$ 
with $p\lef\sq p_0$, 
which \dd{\ls_{\ga_\za}}forces that $p\rig\in t\lef[\uG]$, 
so that $p_0\in\aal{a_\za}$. 
We can assume that $\za=0$ 
(otherwise simply cut all double-names $a_\xi$, $\xi<\za$). 
Then $p\rig\in\Gr^\xi$, and also 
$p_0\in\aal{a_\xi}$ by Lemma~\ref{0ext}\ref{0ext3}, 
for each $\xi<\la$. 

It follows from $p\rig\in\Gr^\xi$ by definition that 
each set $\Gr^\xi$ is 
\dd{\ls_{\ga_\xi}}generic and $\Gl^\xi=t\rig^\xi[\Gr^\xi]$. 
}

To continue with \ref{dq}\ref{dq1}(a), we prove the 
\dd{\ls_\ga}genericity of $\Gr$.

Let $D\sq\ls_\ga$ be a dense subset of $\ls_\ga$, in $\rL$. 
Assume towards the contrary that $\Gr\cap D=\pu$. 
Then there is a condition $p\in \Gl$ 
which \dd{\ls_\ga}forces that $\tl A[\uG]\cap D=\pu$. 
Then $p\in\Gl^\xi$ for some $\xi<\la$, and 
there is a condition $q\in \Gr^\xi$ 
which puts $p$ in $\Gl^\xi=\tr\xi[\Gr^\xi]$  
in the sense that $\ang{q,p}\in\tr{\xi}$. 
As $D$ is dense, there is some $q'\in D$ with $q\sq q'$. 
Then $q'$ belongs to some 
$\ls_{\ga_\eta}$, $\xi<\eta<\la$. 
By Corollary~\ref{qq'}, there is a condition 
$p'\in \ls_{\ga_\eta}$, compatible with $p$ and such 
that $\ang{p',q'}\in \tl\eta$.
Then $p'$ \dd{\ls_\ga}forces $q'\in\tl\eta[\uG]\cap D$, 
while $p$, a compatible condition, forces the opposite, 
which is a contradiction.
\vyk{
Consider an arbitrary \dd{\ls_{\ga_\eta}}generic set 
$H\rig\sq\ls_{\ga_\eta}$ containing $r$, hence, $q$ as well. 
Then, as $q\in\dom{\tr\xi}\sq\dom{\tr\eta}$, the set 
$H\lef=t\rig^\eta[H\rig]$ is \dd{\ls_{\ga_\eta}}generic     
and $H\rig=t\lef^\eta[H\lef]$ by Lemma~\ref{io}. 
Moreover, by Lemma~\ref{0ext},
the sets $H'\rig=H\rig\cap\ls_{\ga_\xi}$ and 
$H'\lef=H\lef\cap\ls_{\ga_\xi}$ are 
\dd{\ls_{\ga_\xi}}generic, and  
$H'\lef=t\rig^\xi[H'\rig]$. 
However $q\in H'\rig$ by construction, hence 
$p\in H'\lef$ by the choice of $q$. 
Thus $p\in H\lef$ as well.

Now consider an \dd{\ls_\ga}generic set 
$K\lef\sq\ls_\ga$ with $H\lef\sq K\lef$; 
then $p\in K\lef$.
The set $K\rig=\tl A[K\lef]$ satisfies 
$H\rig\sq K\rig$, therefore  
$r\in K\rig$ and hence $K\rig\cap D\ne\pu$, 
contrary to the choice of $p$ and $K\lef$. 
This ends the proof of the genericity of $\Gr$ 
and \ref{dq}\ref{dq1}(a).
}%
\vyk{
Finally to establish \ref{dq}\ref{dq3b}, it suffices to 
define a \dd{\ls_\ga}generic set $\Gl\sq \ls_\ga$, 
such that $p\lef\in \Gl$ and $p\rig\in \Gr=t\lef[\Gl]$.
Recall that $p\lef\in\ls_{\ga_0}$.
As ${a}_0\in\trn_{\ga_0}$, there is a 
\dd{\ls_{\ga_0}}generic set $H\lef\sq \ls_{\ga_0}$, 
containing $p\lef$ and such that 
$p\rig\in H\rig=t\lef[H\lef]$. 
Following an argument above, consider a 
\dd{\ls_{\ga}}generic set $\Gl\sq \ls_{\ga}$
with $H\lef=\Gl\cap\ls_{\ga_0}$.
Then $H\rig\sq \Gr=t\lef[\Gl]$, and 
$p\lef\in \Gl$,  $p\rig\in \Gr$ by construction.
}%

\ref{3ext3} Pretty similar argument.
\epf  

\bcor
[in $\rL$]
\lam{resq}
Assume that\/ $c\in\fpc\Om$.
Then 
\ben
\renu
\itla{resq1}
the set\/ $\Xi=\ens{\ga<\Om}{c\res\ga\in\trn_\ga}$ 
is a club in\/ $\Om\,;$ 

\itla{resq2}
if\/ $c$ is full\/ {\rm(Definition~\ref{fn})} then 
$\Xi'=\ens{\ga\in\Xi}{c\res\ga\,\text{\rm\ is full}}$
is a club$;$

\itla{resq3}
if\/ $\Xi''=\ens{\ga\in\Xi}{c\res\ga\,\text{\rm\ is regular}}$
is unbounded in\/ $\Om$ then\/ $\Xi''=\Xi\,.$
\een 
\ecor  
\bpf
\ref{resq1}
That $\Xi$ is closed follows from 
Lemma~\ref{3ext}\ref{3ext1}. 
To prove that $\Xi$ is unbounded,  
let $\al<\Om$ and find a larger ordinal $\ba\in\Xi$.

Recall that to decide a sentence $\Phi$ means to force 
$\Phi$ or to force $\neg\:\Phi$. 

By basic forcing theorems, if $p\in\ls$ then the set 
$$
D_p=\ens{p\in\ls}
{p\text{ decides }q\in\tl c[\uG]\text{ and decides }q\in\tr c[\uG]}
$$
is dense in $\ls$, therefore by the ccc property 
of $\ls$ there is 
an ordinal $\ba$, $\al<\ba<\Om$, such that $D_p$ is dense 
in $\ls_{\ba}$ for all $p\in\ls_\ba$. 
Then $\ba\in\Xi$.

\ref{resq2} easily follows from \ref{resq1}.
To prove \ref{resq3} apply Lemma~\ref{erh}\ref{erh3}. 
\epf

\parf{Superpositions}
\las{sus}

Assume that $\ga\le\Om$ and $a,c\in\trn_\ga$. 
Define 
$$
\bay{rcl}
\tl{a\app c} 
&=& 
\ens{\ang{p',q}\in\ls_\ga\ti\ls_\ga}
{
\sus p\in\ls_\ga\, 
(\ang{p',p}\in \tl c\land \ang{p,q}\in\tl a)},
\\[1.5ex]
\tr{a\app c} 
&=& 
\ens{\ang{q,p'}\in\ls_\ga\ti\ls_\ga}
{
\sus p\in\ls_\ga\, 
(\ang{q,p}\in\tr a\land \ang{p,p'}\in \tr c)}.
\eay
$$
\index{zztaclef@$\tl{a\app c}$}%
\index{zztacrig@$\tr{a\app c}$}%
\index{zzac@$a\app c$}%
\index{superposition!ac@$a\app c$}%
and $a\app c=\ang{\tl{a\app c},\tr{a\app c}}$.

\ble
\lam{sus1}
If\/ $\ga\le\Om$, $a,c\in\trn_\ga$, 
and\/ $G\sq\ls_\ga$, then\/ 
$\tl{a\app c}[G]=\tl a[\tl c[G]]$ and\/ 
$\tr{a\app c}[G]=\tr c[\tr a[G]]$.
\ele
\bpf
Assume that $q\in\tl{a\app c}[G]$. 
Then there is a pair $\ang{p',q}\in \tl{a\app c}$ 
with $p'\in G$. 
By definition there is a condition $p$ such that 
$\ang{p',p}\in \tl c$ and $\ang{p,q}\in\tl a$. 
Then $p\in\tl c[G]$ and hence $q\in \tl a[\tl c[G]]$.
To prove the converse assume that $q\in\tl a[\tl c[G]]$. 
Then there is a pair $\ang{p,q}\in \tl{a}$ 
with $p\in\tl c[G]$, and further 
there is a pair $\ang{p',p}\in \tl{c}$ with $p'\in G$. 
Then $p$ witnesses that $\ang{p',q}\in \tl{a\app c}$, 
and hence $q\in \tl{a\app c}[G]$. 
\epf

\bcor
\lam{sus1c}
Assume that\/ $\ga<\Om$ and\/ $a,b,c\in\ls_\ga$. 
If\/ $a,b$ are equivalent  
(in the sense of Definition~\ref{eqn}) 
then so are\/ $a\app c$ and\/ $b\app c$.\qed 
\ecor

\ble
\lam{sus2}
If\/ $\ga\le\Om$ and\/ $a,c\in\trn_\ga$ 
then the following are equivalent$:$
\bce 
{\rm(1)} $\ran{\tl c}\cap\dom{\tl a}\ne\pu$, \
{\rm(2)} $\ran{\tr a}\cap\dom{\tr c}\ne\pu$, \ 
{\rm(3)} $a\app c\in\trn_\ga$.
\ece
\ele
\bpf
Let $\ran{\tl c}\cap\dom{\tl a}\ne\pu$. 
To prove (3) consider an \dd{\ls_\ga}generic set 
$G'\sq\ls_\ga$, and let $p'\in G'\cap\dom{\tl{a\app c}}$. 
Then $p'\in \dom{\tl c}$, hence $G=\tl b[G']$ is 
an \dd{\ls_\ga}generic set by Lemma~\ref{io}. 
As $p'\in \dom{\tl{a\app b}}$,   
$G\cap\dom{\tl a}\ne\pu$.  
It follows that $H=\tl a[G]$ is 
an \dd{\ls_\ga}generic set.
Finally $H=\tl{a\app c}[G']$ by Lemma~\ref{sus1}.

This argument also proves that $G'=\tr{a\app c}[H]$. 
Thus ${\rm(1)}\imp {\rm(3)}$.

That ${\rm(3)}\imp {\rm(1)}$ is obvious. 
\epf

\bcor
\lam{sus3}
If\/ $\ga\le\Om$, $a,c\in\trn_\ga$, and\/
$c$ is full\/ (in the sense of Definition~\ref{fn}) 
then\/ $a\app c\in\trn_\ga$. 
\ecor
\bpf
By Lemma~\ref{fnd}, 
$\ran{\tl c}=\ran{\tr c}=\ls_\ga$.   
Now use Lemma~\ref{sus2}.   
\epf

Thus if $c\in\trn_\ga$ is a full double-name then 
$a\mto a\app c$ is a map $\trn_\ga\to\trn_\ga$. 
In this case, consider the \rit{inverse} double-name 
\index{double-name!inverse, $c\obr$}%
\index{zzc-1@$c\obr$}%
$c\obr=\ang{\tr c,\tl c}$, let $a\in\trn_\ga$, and 
compare $a$ with $a'=a\app c\app c\obr$. 
On the one hand, we have 
$\tl{a'}[G]=\tl a[\tl c[\tl {c\obr}[G]]]$ for any 
\dd{\ls_\ga}generic set $G$ by Lemma~\ref{sus1}. 
It follows that 
$\tl{a'}[G]=\tl a[\tl c[\tr {c}[G]]]=\tl a[G]$ since 
the successive action of $\tl c$ and $\tr c$ is the 
identity by Lemma~\ref{io}. 
Similarly $\tr{a'}[G]=\tr a[G]$. 
Therefore $a$ and $a'$ are equivalent, and hence 
$\rh a=\rh{a'}$ by Lemma~\ref{rh}, but 
generally speaking we cannot assert that 
straightforwardly $a=a'$. 

To fix this problem, define 
the modified action $a\aqq c=\rh{(a\app c)}$.
\index{zzac*@$a\aqq c$}%
\index{superposition!ac*@$a\aqq c$}%

\ble
\lam{c*c}
Let\/ $\ga<\Om$ and let\/ 
$c\in\trn_\ga$ be a full double-name. 
If\/ $a\in\trn_\ga$ is regular 
(that is, $a=\rh a$) then\/ $b=a\aqq c\in\trn_\ga$, 
$b$ is regular, and\/ $a=b\aqq c\obr$. 
\ele
\bpf
That $b\in\trn_\ga$ follows from Corollary~\ref{sus3}. 
The regularity holds by Lemma~\ref{rh}. 
To prove $a=b\aqq c\obr$, note that both $a$ and 
$b\aqq c\obr$ are regular double-names, and hence it 
suffices, by Lemma~\ref{rh}, 
to prove that $a$ and 
$b\aqq c\obr$ are equivalent. 
However, still by Lemma~\ref{rh}, 
$b\aqq c\obr$ is equivalent to $b\app c\obr$, 
and $b=a\aqq c$ is equivalent to $a\app c$, hence 
overall $b\aqq c\obr$ is equivalent to $a\app c\app c\obr$ 
by Corollary~\ref{sus1c}. 
Finally $a$ is equivalent to $a\app c\app c\obr$, see above. 
\epf


\ble
\lam{cles}
Assume that\/ $\ga<\da\le\Om$,
$c\in\trn_\ga$ and\/ $d\in\trn_\da$ are full double-names, 
$c=d\res\ga$, and\/ $a\in\trn_\ga$, $b\in\trn_\da$. 
Then
\ben
\renu
\itla{cles1} 
if\/ $a\les b$ then\/ $a\app c\les b\app d\,;$  

\itla{cles2} 
if\/ $a,b$ are regular then\/  
$a\les b$ iff\/ $a\aqq c\les b\aqq d$. 
\een
\ele
\bpf
\ref{cles1} is clear since $a\app c$ is monotone on 
both $a$ and $c$. 
As for \ref{cles2}, the implication $\imp$ holds by 
\ref{cles1} and Lemma~\ref{erh} while to prove the 
inverse make use of Lemma~\ref{c*c}.
\epf

\parf{Generic double-names and product forcing} 
\las{gnpf}

By Lemma~\ref{3ext}, we can consider the set 
$\trn=\bigcup_{\ga<\Om}\trn_\ga$ ordered by $\les$ 
as an \dd\Om closed forcing notion in $\rL$ 
(\dd\les bigger double-names are stronger conditions). 
Suppose that $\Ga\sq\trn$ is a \dd\trn generic set 
over $\rL$. 
Then a double-name $\xa=\bigvee\Ga\in\rL[\Ga]$ can be defined 
as in Section~\ref{inc}; we call such double-names 
$\xa=\bigvee\Ga$  
\rit{generic over\/ $\rL$} 
(together with the background generic sets $\Ga$).

Let $\uga$ and $\ua$ be canonical \dd\trn names of 
resp.\ $\Ga$ and $\xa=\bigvee\Ga$.
\index{zzGaund@$\uga$}%
\index{zzAund@$\ua$}%

\bre
\lam{gnar}
As $\rL$ is our default ground model unless otherwise 
specified, the sets $\Ga$ and $\xa=\bigvee\Ga$ 
do not belong to 
$\rL$, however all reals and generally all sets 
$x\sq\ga<\Om$ in $\rL[\Ga]$ belong to $\rL$ 
by Lemma~\ref{3ext}.
It follows that the definition of $\trn_\ga$ ($\ga<\Om$) 
in $\rL$ is absolute for $\rL[\Ga]$. 
That is, if ${a}\in \trn_\ga$ in $\rL$ then it is true 
in $\rL[\Ga]$ that ${a}\in \trn_\ga$. 
And conversely, if ${a}\in \rL[\Ga]$ and it is true 
in $\rL[\Ga]$ that ${a}\in \trn_\ga$ then ${a}\in \rL$ 
and it is true in $\rL$ that ${a}\in \trn_\ga$. 
\ere

\bcor
\lam{gga}
Assume that\/ $\Ga$ is\/ \dd{\trn}generic over\/ $\rL$ 
and\/ $\xa=\bigvee\Ga$. 
Then\/ 
\ben
\renu
\itla{gga1}
it holds in\/ $\rL[\Ga]$ that\/ $\xa$ 
belongs to\/ $\trn_\Om\,;$ 

\itla{gga2}
if\/ $\Gl$ is\/ 
\dd{\ls}generic over\/ $\rL[\Ga]$, 
and\/ $\Gl\cap\dom\tl\xa\ne\pu$, 
then\/ $\Gr$ 
is\/ \dd\ls generic over\/ $\rL[\Ga]$ and\/ 
$\Gl=\tr\xa[\Gr]\,;$ 

\itla{gga3}
if\/ $a\in\trn$, $a\sq\xa$, and\/ $\ga=\abt{a}$ then\/ 
$A\res\ga\in \trn_\ga\cap\Ga$ and\/ $a\les A\res\ga\les A$.
\een
\ecor
\bpf
\ref{gga1}
Remark~\ref{gnar} allows simply to refer to 
Lemma~\ref{3ext}.

\ref{gga2}
Make use of Lemma~\ref{io}.

\ref{gga3}
To prove that $a'=A\res\ga\in \trn_\ga$ and $a\les a'\les A$ 
refer to Lemma~\ref{0ext}\ref{0ext3}.
To prove that $a'\in \Ga$ note that by Lemma~\ref{3ext} 
there is some $c\in\Ga$ which decides each $b\in\ls_\ga$ 
to belong or not to belong to $\Ga$; then $a'\sq c$.
\epf

\vyk{
The corollary obviously refers to the product forcing 
$\ls\ti\trn\in\rL$. 
Yet it will be more convenient to consider a certain 
subforcing of $\ls\ti\trn$. 

\bdf
\lam{lsp:d}
A double-name $\ang{p,a}\in\ls\ti\trn$ is \rit{compatible} 
if $p\in\dom\tl a$. 
Let $\lso$ be the set of all compatible double-names. 
\edf

\ble
\lam{ecom}
If\/ $p\in\ls$ then there is a double-name\/ 
$\ang{p,a}\in\lso$.
\ele
\bpf
Choose an ordinal $\ga<\Om$ such that $p\in\ls_\ga$. 
Define $a=\ang{t,t}$, where 
$t=\ens{\ang{q,q}}{q\in\ls_\ga}$.
\epf

The next lemma shows that the compatibility 
restriction is rather innocuous.  

\bcor
\lam{ggb}
Let\/ $\Gl\ti\Ga$ be\/ \dd{(\lsp)}generic over\/ $\rL$, 
and\/ $\xa=\bigvee\Gl$. 
Then the set\/ $\Gl\ti\Ga$ contains a double-name in\/ $\lso$ 
if and only if\/  
the set\/ $\Gr=\tl\xa[\Gl]$ contains\/ $\pr\xa$.

iff\/ $\Gl\cap\dom\tl\xa\ne\pu$.
\ecor
\bpf
To prove $\imp$, assume that $\ang{p,a}\in\lso$ and 
$p\in\Gl$, $a\in\Ga$. 
By definition 
$\pl\xa=\pl a$ and $\pr \xa=\pr a$. 
Let $\ga<\Om$ be such that $a\in\trn_\ga$.
Then $H=\Gl\cap\ls_\ga$ is \dd{\ls_\ga}generic    
and $\tl a[H]\sq\tl\xa[\Gl]=\Gr$ by 
Lemma~\ref{0ext}\ref{0ext1}. 
On the other hand, there is $p\in\Gl$ such that 
$\ang{p,a}\in\lsp$, hence $p\in H$ and $p\in\aal {H}$. 
It follows that $\pr \xa=\pr a\in \tl a[H]$, as required.

To prove $\mpi$, assume to the contrary that a pair 
$\ang{p,a}\in\Gl\ti\Ga$ forces that $\uG\ti\uga$ contains no 
pairs in $\lso$ but $\pr\ua\in\tl\ua[\uG]$. 
Let $\ga<\Om$ be the least satisfying $a\in\trn_\ga$. 
We can assume that $p\in\ls_\ga$, as otherwise $a$ 
can be extended by Lemma~\ref{ext1}.
Note that $p$ \dd{\ls_\ga}forces that 
$\pr a\nin\tl a[\uG]$. 
(As otherwise there is a stronger 
$q\in\ls_\ga$ which forces $\pr a\in\tl a[\uG]$, hence 
$q\in\aal a$ and $\ang{q,a}\in\lso$, contrary to the 
choice of $p$.) 
The set $H=\Gl\cap\ls_\ga$ is \dd{\ls_\ga}generic as 
above, and contains $p$, and hence 
$\pr a\nin\tl a[H]$.
However $\pr a=\pr\xa$ and 
$\tl a[H]=\tl\xa[\Gl]\cap\ls_\ga$
by Lemma~\ref{0ext}, so that $\pr\xa\nin\tl\xa[\Gl]$, 
which contradicts the contrary assumption.
\epf

}

\parf{The first ingredient} 
\las{gna}

Generic double-names and forcing with $\lsp$ 
enable us to carry out 
the first main step towards 
Theorem~\ref{mt}. 

In $\rL$, let $\ho$ be the set of all sets $x$ such 
that the transitive closure $\tc(x)$ has cardinality 
$\card{\tc(x)}<\Om$ strictly.

\bbla
\lam{ass1}
Thus suppose that $G_0\sq\ls$ is a \dd\ls generic 
\index{zzG0@$G_0$}
set over $\rL$, let $\cX\in\rL[G_0]$, 
\index{zzXc@$\cX$}
and it is true in $\rL[G_0]$ that $\cX$ is
a countable \od\ non-empty set of sets of reals. 
%
There is a formula $\vpi(\cdot,\pi)$ with some 
\index{zzphi@$\vpi(\cdot,\pi)$}
$\pi\in\Ord$ as the only parameter, such that it 
is true in $\rL[G_0]$ that $\cX$ is the only set $x$
satisfying $\vpi(x,\pi)$.

There is a sequence
$\dr=\sis{U_n}{n\in\om}\in\rL$ of names $U_n\in\rL$, 
\index{zzud@$\du$}%
\index{zzUn@$U_n$}%
such that
$\cX=\bint\bbc{G_0}:=\ens{\bint{U_n}{G_0}}{n\in\om}$.
Each $U_n$ can be assumed to be 
an \dd\ls name of a set of reals, 
that is, in $\rL$, $U_n\sq\ls\ti \mathbb T$, where $T$ 
is the set of all \dd\ls names for reals. 
Furthermore, according to the \dd\Om cc property of the 
forcing $\ls$, 
each \dd\ls name for a real can be assumed to be 
a set in $\ho$. 
Therefore we shall wlog assume that $U_n\sq\ho$ for all $n$.

Anyway there is a condition\/ $\bp\in G_0$ which \dd\ls 
\index{zzp-@$\bp$}%
forces over $\rL$ that  
``$\bint\bbc{\uG}$ is the only set $x$ 
satisfying\/ $\vpi(x,\pi)$, and 
$\bint\bbc{\uG}$ is a set of sets of reals''.
Let $\bga<\Om$ be the least ordinal satisfying 
\index{zzga-@$\bga$}%
$\bp\in\ls_{\bga}$.
\qed
\ebla

Let a \dd\bp\rit{pair} be any pair $\ang{p,a}\in\lsp$ 
\index{zzp-pair@$\bp$-pair}%
\index{pair!zzp-pair@$\bp$-pair}%
such that $\bp\sq p\in \dom\tl a$ and $p$ 
\dd{\ls_\ga}forces that $\bp\in\tl a[\uG]$, where 
$\ga=\abt{a}$. 

\bre
\lam{ppe}
Let  
$\baa=\ide\bga$. 
Then $\ang{\bp,\baa}$ is a \dd\bp pair; 
$\bp$ \dd{\ls_\bga}forces that $\tl \baa[\uG]=\uG$. 
\ere

\ble
\lam{lsp}
Let\/ $\ang{p,a}\in\lsp$ be a\/ \dd\bp pair, 
$q\in\ls$, $b\in\trn$, $p\sq q$, $a\les b$. 
There is a double-name\/ $c\in\trn$ such that\/ $b\les c$ 
and\/ $\ang{q,c}$ is a\/ \dd\bp pair.
\ele
\bpf
If $q\in\ls_\ga$, where $\ga=\abt{b}$, 
then to define $c$ add to $\tl b$ 
all pairs $\ang{q,r}$ such 
that already $\ang{p,r}\in b$.
We claim that $\ang{q,c}$ is a\/ \dd\bp pair. 
Indeed if $\Gl\sq \ls_\ga$ is generic then easily 
(*) $\tl c[\Gl]=\tl b[\Gl]$, hence $c\in\trn_\ga$. 
Further $\bp\sq p\sq q\in\dom \tl c$ by construction.
Finally $q$ \dd{\ls_\ga}forces that $\bp\in\tl a[\uG]$
because so does $p$, and we can replace $\tl a$ by $\tl c$ 
since $a\sq b\sq c$.

If $q\nin\ls_{\ga}$ then still 
$q\in\ls_\da$ for some $\da$, $\ga<\da<\Om$. 
Use Lemma~\ref{ext1} to get a double-name $b'\in\trn_\da$
with $b\les b'$, and argue as in the first case.
\epf 

\bte
\lam{1mt}
Suppose that\/ $\Gl\ti\Ga$ is a\/ \dd{\lsp}generic set 
over\/ $\rL$, $\xa=\bigvee\Ga$, 
and\/ $\ang{p,a}\in \Gl\ti\Ga$ is a\/ \dd\bp pair. 
Then\/ 
\ben
\renu
\itla{1mt1}
$p,\bp\in\Gl$, $\bp\in\Gr=\tl\xa[\Gl]$, and\/ 
$\Gr$ is\/ \dd\ls generic over\/ $\rL[\Ga]\,;$ 
 
\itla{1mt2}
$\bint\bbc\Gl=\bint\bbc\Gr$ --- 
in other words, any\/ \dd\bp pair\/ $\ang{p,a}$
\dd{(\lsp)}forces\/ 
$\bint\bbc\uG=\bint\bbc{\tl\ua[\uG]}$  
over\/ $\rL$.   
\een
\ete
\bpf
\ref{1mt1} 
To prove the genericity apply Corollary~\ref{gga}. 

To prove \ref{1mt2} suppose otherwise.  
Then there is a pair $\ang{q,b}$ in $\lsp$ 
with $p\sq q$, $a\les b$,
which \dd{(\lsp)}forces 
$\bint\bbc\uG\ne\bint\bbc{\tl\ua[\uG]}$, 
that is   
\ben
\fenu
\atc
\itla{rl}
if $\Gl\ti\Ga$ is a \dd{(\lsp)}generic set 
over\/ $\rL$ containing $\ang{q,b}$, 
$\xa=\bigvee\Ga$, and $\Gr=\tl\xa[\Gl]$, 
then $\bint\bbc\Gl\ne\bint\bbc\Gr$. 
\een

Let $\cL\in\rL$ be an elementary
submodel of a large model, such that $\ho\sq\cL$,
$\Om$ and $\pi$ belong to $\cL$, $\card{\cL}=\Om$ in $\rL$,
and $\cL$ is an elementary submodel of $\rL$ \vrt\
all $\is{}{100}$  formulas.
Let $\cL'\in\rL$ be the Mostowski collapse of $\cL$; 
still $\card{\cL'}=\Om$ in $\rL$.
Note that $\cL'$ is a transitive model of Zermelo with choice,
and the collapse map $\phi:\cL\onto\cL'$ is the identity 
on $\ho$, hence even on $\pws\ho\cap\cL$. 
In particular, 
$\phi(\Om)=\Om$, 
$\phi(\bbc)=\bbc$, $\phi(U_n)=U_n$ for all $n$, 
$\phi(\ls)=\ls$, $\phi(\trn)=\trn$, $\ho\sq\cL'$, and even 
$\pws\ho\cap\cL\sq\cL'$. 


By the elementary submodel property, 
$\ang{q,b}$ 
still \dd{(\lsp)}forces over $\cL'$ that 
$\bint\bbc\uG\ne\bint\bbc{\tl\ua[\uG]}$  
--- that is
\ben
\fenu
\atc
\atc
\itla{cl}
if $\Gl\ti\Ga$ is a \dd{(\lsp)}generic set 
over\/ $\cL'$ containing $\ang{q,b}$, 
$\xa=\bigvee\Ga$, 
and $\Gr=\tl\xa[\Gl]$, 
then $\bint\bbc\Gl\ne\bint\bbc\Gr$. 
\een
To infer a contradiction, note that since $\card{\cL'}=\Om$ 
in $\rL$, by Lemma~\ref{3ext} there exists a set $\Ga\in\rL$, 
\dd{\trn}generic over\/ $\cL'$ and containing $b$, hence 
containing $a$ as well. 
We underline that $\Ga\in\rL$, and then $\xa=\bigvee\Ga$ 
belongs to $\rL$, too.
Let $\Gl\sq\ls$ be a set \dd{\ls}generic over $\rL$, 
hence over $\cL'[\Ga]$ as well, 
and containing $q$, and then containing $p$. 
Then the set $\Gr=\tl\xa[\Gl]$ is \dd{\ls}generic 
over $\rL$ and over $\cL'[\Ga]$ by Lemma~\ref{io}, and 
in addition, $\bint\bbc\Gl\ne\bint\bbc\Gr$ by \ref{cl}.

Recall that $\ang{p,a}$ also belongs to $\Gl\ti\xa$. 
Therefore $\bp\in\Gl\cap\Gr$ by \ref{1mt1}. 
Thus $\Gl$ and $\Gr$ are \dd\ls generic sets over $\rL$ 
and both contain $\bp$, 
$\nt\bbc{\Gl}$ is the only set $x$ satisfying 
$\vpi(x,\pi)$ in $\rL[\Gl]$ 
while $\nt\bbc{\Gr}$ is the only set $x$ satisfying 
$\vpi(x,\pi)$ in $\rL[\Gr]$. 
However $\rL[\Gl]=\rL[\Gr]$ 
(because $\Gr=\tl\xa[\Gl]$, $\Gl=\tr\xa[\Gr]$, 
and $\xa\in\rL$), 
while on the other hand  
$\bint\bbc\Gl\ne\bint\bbc\Gr$, 
which is a contradiction.
\epf

\parf{Stabilizing pairs and second ingredient}
\las{simt}

Let a \rit{stabilizing \dd\bp pair} be any \dd\bp pair 
\index{zzp-pair@$\bp$-pair!stabilizing}%
\index{pair!zzp-pair@$\bp$-pair!stabilizing}%
$\ang{\hp,\ha}\in\lsp$ 
which, for some\/ $n$, \dd{(\lsp)}forces\/ 
$\bint{U_0}\uG =\bint{U_{n}}{\tl\ua[\uG]}$ over\/ $\rL$.

\bcor
\lam{mt1c}
If\/ $\Gl$ is an\/ \dd{\ls}generic set over\/ $\rL$ 
containing\/ $\bp$, 
then there is a stabilizing\/ \dd\bp pair\/ 
$\ang{\hp,\ha}\in\lsp$ with\/ $\hp\in\Gl$. 
\ecor
\bpf
Let $\baa=\ide\bga$, see Remark~\ref{ppe}. 
Let $\Ga\sq\trn$ be a set \dd\trn generic over 
$\rL[\Gl]$ and containing $\baa$, so that $\Gl\ti\Ga$ 
is \dd{(\lsp)}generic. 
Let $\xa=\bigvee\Ga$. 
Then the set $\Gr=\tl\xa[\Gl]$ satisfies  
$\bint\bbc\Gl=\bint\bbc\Gr$
by Theorem~\ref{1mt}.
Therefore there is a number $n\in\om$ such that 
$\bint{U_0}\Gl=\bint{U_{n}}\Gr$. 
Then there is a 
stronger pair $\ang{\hp,\ha}\in\Gl\ti\Ga$
($\bp\sq\hp$ and $\baa\les \ha$) 
which \dd{(\lsp)}forces\/ 
$\bint{U_0}\uG=\bint{U_{n}}{\tl\ua[\uG]}$.
We can assume that $\ang{\hp,\ha}$ is a \dd\bp pair, 
by Lemma~\ref{lsp}.
\epf

\bpro
\lam{mt1d}
Let $\ang{\hp,\ha}\in\lsp$ be a stabilizing\/ \dd\bp pair. 
Assume that\/ $\Gl\ti\Ga$, $\Gl'\ti\Ga'$ are sets\/ 
\dd{(\lsp)}generic over\/ $\rL$ 
and containing\/ $\ang{\hp,\ha}$, 
$\xa=\bigvee\Ga$, $\xa'=\bigvee{\Ga'}$, and\/ 
$\tl\xa[\Gl]=\tl{\xa'}[\Gl']$.
Then\/ $\bint{U_0}\Gl=\bint{U_0}{\Gl'}$. 
\epro
\bpf
By definition,
$\bint{U_0}\Gl =\bint{U_{n}}{\tl\xa[\Gl]}$ and 
$\bint{U_0}{\Gl'} =\bint{U_{n}}{\tl{\xa'}[\Gl']}$
for one and the same $n$.
\epf

The second ingredient in the proof of Theorem~\ref{mt}
will be the following:

\bte
\lam{mt2}
Assume that\/ $\ang{\hp,\ha}\in\lsp$ is a  
stabilizing\/ \dd\bp pair, 
$\hga<\Om$, $\ha\in\trn_\hga$, $\hp\in\ls_\hga$, 
$\Gl,\Gl'\sq\ls$ are\/ \dd\ls generic sets 
over\/ $\rL$ containing\/ $\hp$, 
$\Gl\cap\ls_\hga=\Gl'\cap\ls_\hga$, 
and\/ $\rL[\Gl]=\rL[\Gl']$. 
Then\/ $\bint{U_0}\Gl=\bint{U_0}{\Gl'}$.
\ete

Let's show how this implies Theorem~\ref{mt}. 
The proof of Theorem~\ref{mt2} itself will follow in 
the next sections. 

\bpf[Theorem~\ref{mt} from Theorem~\ref{mt2}]
We argue in the assumptions and notation of \ref{ass1}. 
Let $\Gl=G_0$, so that $\bp\in\Gl$ by \ref{ass1}. 
Then by Corollary~\ref{mt1c}, 
there is a stabilizing \dd\bp pair 
$\ang{\hp,\ha}\in\lsp$ such that $\hp\in\Gl$. 
Pick $\hga<\Om$ such that $\ha\in\trn_\hga$ and 
$\hp\in\ls_\hga$. 
Consider, in $\rL[\Gl]$, the set $\cG$ of all sets 
$G\sq\ls$, \dd\ls generic over $\rL$ and satisfying 
$\rL[G]=\rL[\Gl]$, $\hp\in G$, and 
$G\cap\ls_\hga=\Gl\cap\ls_\hga$.
In particular $\Gl\in\cG$.
The only essential parameter of the definition of $\cG$ 
which is not immediately \od\ --- is $\Gl\cap\ls_\hga$.
However $\Gl\cap\ls_\hga$ itself, as basically any subset of 
any $\ls_\ga$, $\ga<\Om$, is \rod\ in the Solovay model. 
We conclude that $\cG$ is \rod\ in $\rL[\Gl]$.

On the other hand, suppose that $G\in\cG$.
Then $\bint{U_0}\Gl=\bint{U_0}{G}$ 
by Theorem~\ref{mt2}.
Therefore the set $\bint{U_0}\Gl$ can be defined as
$\bint{U_0}G$ for some / every $G\in\cG$.
This witnesses that $\bint{U_0}\Gl$ is \rod\ in $\rL[\Gl]$,
because so is $\cG$ by the above.
Thus the set $\cX=\bint\bbc\Gl$ contains a \rod\ element.
It follows that $\cX$ contains an \od\ element, by
Lemma~\ref{mtl}, as required.\vom

\epF{Thm~\ref{mt} mod Thm~\ref{mt2}}


\parf{Final} 
\las{mod}

Here we prove Theorem~\ref{mt2} 
and finally prove Theorem~\ref{mt}.
\index{zzzz@{\ubf separate index for Section~\ref{mod}}|see{below}}
%
{\ubf We argue in the assumptions and notation of 
Theorem~\ref{mt2}.}
That is, 
\ben
\nenu
\itla{mod1}
$\ang{\hp,\ha}\in\lsp$ is a stabilizing \dd\bp pair, 
\imar{mod1}
$\hga<\Om$, $\ha\in\trn_\hga$, $\hp\in\ls_\hga$, 
\index{zzpha@$\hp$}%
\index{zzaha@$\ha$}%
\index{zzgaha@$\hga$}%
\index{zzzzpha@$\hp$}%
\index{zzzzaha@$\ha$}%
\index{zzzzgaha@$\hga$}%
the sets $\Gl,\Gl'\sq\ls$ are\/ \dd\ls generic   
\index{zzGlef@$\Gl$}%
\index{zzGlefp@$\Gl'$}%
\index{zzzzGlef@$\Gl$}%
\index{zzzzGlefp@$\Gl'$}%
over\/ $\rL$ and both contain $\hp$, and in addition
$\Gl\cap\ls_\hga=\Gl'\cap\ls_\hga$, 
$\rL[\Gl]=\rL[\Gl']$.
\een 
In this assumption, we have to prove that 
$\bint{U_0}\Gl=\bint{U_0}{\Gl'}$.
Working towards this goal, our plan will be to find:
\ben
\fenu
\itla{*}
sets $\Ga,\Ga'\sq\trn$, 
\imar{*}
\dd\trn generic over $\rL[\Gl]=\rL[\Gl']$, 
containing $\ha$, 
and satisfying $\tl\xa[\Gl]=\tl{\xa'}[\Gl']$, where 
$\xa=\bigvee\Ga$ and $\xa'=\bigvee\Ga'$;
\een
then the products $\Gl\ti\Ga$ and $\Gl'\ti\Ga'$ will be  
\dd{(\lsp)}generic over $\rL$ 
and containing $\ang{\hp,\ha}$,
so that $\bint{U_0}\Gl=\bint{U_0}{\Gl'}$ follows by 
Proposition~\ref{mt1d}, accomplishing the proof 
of Theorem~\ref{mt2}.

By Theorem~\ref{GG} there is a double-name\/ 
$C\in\trn_\Om$ in $\rL$, such that\/ 
\index{zzC@$C$}%
\index{double-name!C@$C$}%
\ben
\nenu
\atc
\itla{mod2}
$C$ is full, 
$\tl C=\tr C$, 
\imar{mod2}%
$\Gl=\tl C[\Gl']$, 
and $\Gl'=\tr{C}[\Gl]$. 
\vyk{ --- 
and we can assume that each condition 
$p\in\dom{\tl C}$ \dd\ls forces that $\hp\in\tl C[\uG]$,  
while each condition  
$q\in\dom{\tr C}$ \dd\ls forces that $\hp\in\tr C[\uG]$.
}%
\een

As $\Gl\cap\ls_\hga=\Gl'\cap\ls_\hga$, we 
can further assume that 
\ben
\nenu
\atc
\atc
\itla{mod3}
the restricted double-name $C\res\hga$ 
\imar{mod3}
coincides with $\ide\hga$ of Example~\ref{idg}, 
so that $C\res\hga\in\ls_\hga$ is full and regular, 
and $\tl{C\res\hga}[G]=\tr{C\res\hga}[G]=G$ for all $G$.
\een

Let $\Ga$ be any set $\Ga\sq\trn$ with $\ha\in\Ga$, 
\dd\trn generic over $\rL[\Gl]$.
\index{zzGa@$\Ga$}%
\index{zzzzGa@$\Ga$}%
Then 
$\xa=\bigvee\Ga\in\trn_\Om$ 
in $\rL[\Ga]$ by Corollary~\ref{gga},  
\index{zzA@$\xa$}%
\index{zzzzA@$\xa$}%
and $\bp\sq\hp\in \dom\xa$ since $\ha\in\Ga$. 

\bcor
\label{resr}
\ben
\renu
\itla{resr1}
The set\/ 
\imar{resr}
$X=\ens{\ga<\Om}{A\res\ga\in\trn_\ga}\in\rL[\Ga]$ 
\index{zzXYZ@$X,Y,Z$}%
\index{zzzzXYZ@$X,Y,Z$}%
is a club in\/ $\Om$, and if\/ $\ga\in X$ then\/ 
$A\res\ga$ is regular$;$  

\itla{resr2}
the set\/ 
\imar{resr}
$Y=\ens{\ga<\Om}{C\res\ga\in\trn_\ga\,
\text{\rm\ and\/ $C\res\ga$ is full}}\in\rL
$ 
is a club in\/ $\Om\,;$  

\itla{resr3}
therefore\/ $Z=\ens{\ga\in X\cap Y}{\hga\le\ga}$ is a club, and 
in addition\/ $\hga\in Z\,.$  
\een 
\ecor
\bpf
To prove \ref{resr1} and \ref{resr2} apply 
Corollary~\ref{resq}; the unboundedness condition in 
\ref{resq}\ref{resq3} follows from the genericity of $\Ga$ 
and the density of the set of all regular double-names 
$a\in\trn$ by Lemma~\ref{rh}\ref{rh2}.

Claim $\hga\in Z$ in \ref{resr3} follows from \ref{mod3}.
%
\epf

Now suppose that $\ga\in Y$, hence $C\res\ga\in\trn_\ga$ 
is full. 
Let $a\in\trn_\ga$ be regular. 
Define $a\aqq C=a\aqq(C\res\ga)$ (see Section~\ref{sus}).

\ble
\lam{1gen}
The map\/ $a\mto a\aqq C$ is a\/ \dd\les preserving 
bijection of the set\/
$\try=
\ens{a\in \trn}{a\,\text{ is regular}\land\abt{a}\in Y}$ 
onto itself, satisfying\/ $a\aqq C\aqq C=a$.%
\index{zzDNYreg@$\try$}%
\index{zzzzDNYreg@$\try$}%
\ele
\bpf
If $a\in \try$ and $\ga=\abt a$ then 
$a\aqq C=a\aqq(C\res\ga)$ belongs to $\trn_\ga$ and 
is regular by Lemma~\ref{c*c}, 
hence $a\aqq C\in \try$. 
If $\da>\ga$ is a bigger ordinal 
\pagebreak[3] 
still in $Y$, and 
$b\in \try$, $\da=\abt b$, then 
$a\les b$ iff $a\aqq C\les b\aqq C$ 
by Lemma~\ref{cles}\ref{cles2}.
Finally $a\aqq C\aqq C=a$ 
holds still by Lemma~\ref{c*c}, because 
$C\obr = C$ (that is, $\tl C=\tr C$) by \ref{mod2}.
\epf

In particular, if $\ga\in Z$ then 
$A\res\ga\in \try$, and hence 
$(A\res\ga)\aqq C\in\try$ is 
a regular double-name. 
Thus $\sis{(A\res\ga)\aqq C}{\ga\in Z}\in\rL[\Ga]$ 
is a \dd\les increasing sequence of regular double-names. 
The following is a key fact.

\ble
\lam{2gen}
The sequence\/ $\sis{(A\res\ga)\aqq C}{\ga\in Z}$ 
is \dd\trn generic over\/ $\rL[\Gl]=\rL[\Gl']$, 
in the sense that if a set\/ $D'\sq\trn$, 
$D'\in\rL[\Gl]$, is open dense in\/ $\trn$ then there is an 
ordinal\/ $\ga\in Z$ such that\/ $(A\res\ga)\aqq C\in D'$.
\ele
\bpf
The set 
$\Da'=D'\cap\try$ belongs to $\rL[\Gl]$ and still 
is dense in $\trn$ by Lemma~\ref{rh}\ref{rh2}.
Therefore its \dd Cimage $\Da=\ens{a\aqq C}{a\in\Da'}$
still belongs to $\rL[\Gl]$ and is dense in $\trn$ by 
Lemma~\ref{1gen}.
It follows by the genericity of $\Ga$ that 
$A\res\ga\in\Da$ for some $\ga\in Z$. 
Then $a=(A\res\ga)\aqq C\in\Da'$, 
since $a\aqq C=A\res\ga$ by Lemma~\ref{1gen}.
\epf

\bcor
\lam{3gen}  
The set\/ 
$\Ga'=\ens{a\in\trn}{\sus\ga\in Z\,(a\les (A\res\ga)\aqq C)}$ 
\index{zzGap@$\Ga'$}%
\index{zzzzGap@$\Ga'$}%
is\/ \dd\trn generic over\/ $\rL[\Gl]=\rL[\Gl']$.\qed
\ecor

Let us check the other intended properties of $\Ga'$ as in 
\ref{*}.

To see that $\ha\in\Ga'$, recall that 
$\ha\in\Ga\cap\trn_\hga$. 
It follows by Corollary~\ref{gga}\ref{gga3} that 
$\ha\les a=\xa\res\hga$.
However $\hga\in Z$ by Corollary~\ref{resr}\ref{resr3}. 
We conclude that $\ha\aqq C\in \Ga'$. 
Finally $\ha\aqq C=\ha\aqq (C\res\hga)=\ha$ since 
$C\res\hga=\ide\hga$  by \ref{mod3}. 
Thus $\ha\in\Ga'$, as required.

Finally prove that $\tl\xa[\Gl]=\tl{\xa'}[\Gl']$, 
where $\xa=\bigvee\Ga$ and $\xa'=\bigvee\Ga'$.
\index{zzAp@$\xa'$}%
\index{zzzzAp@$\xa'$}%
It suffices to show that if $\ga\in Z$ then 
$$
\tl{\xa\res\ga}[\Gl\cap\ls_\ga]=
\tl{\xa'\res\ga}[\Gl'\cap\ls_\ga]\,.
\eqno(5)
$$ 
However by construction 
$\xa'\res\ga=(\xa\res\ga)\aqq C=(\xa\res\ga)\aqq (C\res\ga)$, 
and on the other hand 
$\tl{(\xa\res\ga)\aqq (C\res\ga)}[G]=
\tl{\xa\res\ga}[\tl{C\res\ga}[G]]$ for all $G$ by 
Lemma~\ref{sus1},  
therefore (5) is equivalent to 
$$
\tl{\xa\res\ga}[\Gl\cap\ls_\ga]=
\tl{\xa\res\ga}[\tl{C\res\ga}[\Gl'\cap\ls_\ga]]\,,
$$ 
which obviously follows from
$$
\Gl\cap\ls_\ga=
\tl{C\res\ga}[\Gl'\cap\ls_\ga]\,,
$$ 
and this is a corollary of the equality $\Gl=\tl{C}[\Gl']$ 
in \ref{mod2} by Lemma~\ref{0ext}\ref{0ext1}(b).\vom

\qeD{Theorem~\ref{mt2}}\vom 

This also completes the proof of Theorem~\ref{mt} 
(see the end of Section~\ref{simt}).\vom                                  

\qeD{Theorem~\ref{mt}}


\parf{Conclusive remarks} 
\las{konk}

\bvo
\lam{voo}
Is Theorem~\ref{mt} true for {\bfit arbitrary} sets $\cX$, 
not necessarily sets of reals?
In this general case, the proof given above fails 
in the proof of Theorem~\ref{1mt}, since it is not true 
anymore that $U_n\sq\ho$ and $\phi(U_n)=U_n$.
\evo

It follows from Theorem~\ref{mt} that, in the Solovay model, 
any \od\ set $\cX$ of sets of reals 
containing non-\od\ elements is \rit{uncountable}. 
If moreover $\cX$ is a set {\ubf of reals} then in fact 
$\cX$ contains a perfect subset and hence has cardinality 
$\mathfrak c$ by a profound theorem in \cite{solo}.
Does this stronger result reasonably generalize to  
sets of sets of reals and more complex sets?

\bco
\lam{vo}
It is true in the Solovay model that if $\cX$ is 
an \od\ set then 
\ben
\Renu
\itla{vo1}
if $\cX$ contains only \od\ elements then it is 
\dd\od wellorderable;

\itla{vo2}
if $\cX$ contains only \rod\ elements, among them 
at leat one non-\od\ element, then $\cX$ includes 
a \rod-image of the continuum $\dn\,;$ 

\itla{vo3}
if $\cX$ contains a non-\rod\ element then $\cX$ 
has cardinality $\ge 2^{\mathfrak c}\,$. 
\een
The set of all \dd\ls generic sets over $\rL$ is 
a less trivial example of a set of type \ref{vo3} in the 
Solovay model. 
\eco

A proof of \ref{vo3} would be an alternative (and perhaps 
simpler) proof of Theorem~\ref{mt} of this paper. 

It remains to note that 
Caicedo and Ketchersid \cite{cai} 
obtained a somewhat similar trichotomy result in 
in a strong determinacy assumption.



\bibliographystyle{plain}
{\small
\bibliography{u}
\printindex
}

\end{document}